\documentclass[12pt] {article}

\usepackage{amssymb}
\usepackage{color}

\begin{document}
\newtheorem{Theorem}{Theorem}[section]
\newtheorem{Proposition}[Theorem]{Proposition}
\newtheorem{Lemma}[Theorem]{Lemma}
\newtheorem{Example}[Theorem]{Example}
\newtheorem{Corollary}[Theorem]{Corollary}
\newtheorem{Fact}[Theorem]{Fact}
\newtheorem{Conjecture}[Theorem]{Conjecture}
\newenvironment{Definition} {\refstepcounter{Theorem} \medskip\noindent
 {\bf Definition \arabic{section}.\arabic{Theorem}}\ }{\hfill}
\newenvironment{Remarks} {\refstepcounter{Theorem}
\medskip\noindent {\bf Remarks
\arabic{section}.\arabic{Theorem}}\ }{\hfill}

\newenvironment{Exercise} {\medskip\refstepcounter{Theorem}
     \noindent {\bf Exercise
\arabic{section}.\arabic{Theorem}}\ } {\hfill}

\newenvironment{Proof}{{\noindent \bf Proof\ }}{\hfill$\Box$}

\newenvironment{claim} {{\smallskip\noindent \bf Claim\ }}{\hfill}

\def \blue {\color{blue}}
\def\red{\color{red}}
\def \L {{\cal L}}

\def \B {{\cal B}}
\def \mod {{\rm mod \ }}
\def \iso {\cong}
\def \Lor  {\L_{\rm or}}
\def \Lr {\L_{\rm r}}
\def \Lg {\L_{\rm g}}
\def \I {{\cal I}}
\def \M {{\cal M}}\def\N {{\cal N}}
\def \E {{\cal E}}
\def \Proj {{\mathbb P}}
\def \H {{\cal H}}
\def \x {\times}
\def \Stab {{\rm Stab}}
\def \Z {{\mathbb Z}}
\def \V {{\mathbb V}}
\def \C {{\mathbb C}} \def \Cexp {\C_{\rm exp}}
\def \R {{\mathbb R}}
\def \Q {{\mathbb Q}}\def \K {{\mathbb K}}
\def \F {{\cal F}}\def \A {{\mathbb A}}
 \def \X {{\mathbb X}}
\def \G {{\mathbb G}}
\def\HH {{\mathbb H}}
\def \Nn {{\mathbb N}}\def \Nn {{\mathbb N}}
\def\D {{\mathbb D}}
\def \hat {\widehat}
\def \bar{\overline}
\def \Spec {{\rm Spec}}
\def \bul {$\bullet$\ }
\def\proves {\vdash}
\def \Co {{\cal C}}
\def \ACFp {{\rm ACF}_p}
\def \ACF0 {{\rm ACF}_0}
\def \ee {\prec}
\def \Diag {{\rm Diag}}
\def \Diage {{\rm Diag}_{\rm el}}
\def \DLO {{\rm DLO}}
\def \d {{\rm depth}}
 \def \dist {{\rm dist}}
\def \P {{\cal P}}
\def \ds {\displaystyle}
\def \Fp {{\mathbb F}_p}
\def \acl {{\rm acl}}
\def \dcl {{\rm dcl}}

\def \dom {{\rm dom}}
\def \tp {{\rm tp}}
\def \stp {{\rm stp}}
\def \Th  {{\rm Th}}
\def\< {\Lngle}
\def \> {\rangle}
\def \n {\noindent}
\def \minusdot{\hbox{\ {$-$} \kern -.86em\raise .2em \hbox{$\cdot \
$}}}
\def\exp {{\rm exp}}\def\ex {{\rm ex}}
\def \td {{\rm td}\ }
\def \ld {{\rm ld}}
\def \span {{\rm span}}
\def \tilde {\widetilde}
\def \d {\partial}
\def \del {\partial}
\def \cl {{\rm cl}}
\def \acl {{\rm acl}}
\def \cN {{\cal N}}
\def \Qalg {{\Q^{\rm alg}}}
\def \th {^{\rm th}}
\def \deg { {\rm deg} }
\def\hat {\widehat}
\def\li {\L_{\infty,\omega}}
\def\lo {\L_{\omega_1,\omega}}
\def\lk {\L_{\kappa,\omega}}
\def \ee {\prec}\def \bSigma {{\mathbf\Sigma}}
\def \mod {\ {\rm mod\ }}
\def \Tor {{\rm Tor}}
\def \td {{\rm td}}
\def \dim {{\rm dim\ }}
 \def \| {\kern -.3em \restriction \kern -.3em}
 \def \lc{\lceil}
 \def \rc{\rceil}
 \def\rcl{{\rm rcl}}
 \def \Pr {{\rm Pr}}\def \PA {{\rm PA}}
 \def \^ {\hat{\  \ }}
 \def \red {\color{red}}
\def \blue {\color{blue}}

\title{Representing Scott Sets in Algebraic Settings\thanks{This work was begun at a workshop on computable stability theory held at the American Institute of Mathematics in August 2013.}}

\author{Alf Dolich\\ Kingsborough Community College\and Julia F. Knight\\University of Notre Dame\and Karen Lange \thanks {Partially supported by National Science Foundation grant DMS-1100604}\\ Wellesley College\and David Marker\\University of Illinois at Chicago
}
\date{}
\maketitle

\section{Introduction}    
    
Recall that $S\subseteq 2^\omega$ is called a {\em Scott set} if and only if:

i) $S$ is a Turing ideal, i.e., if $x,y\in S$ and $z\le_T x\oplus y$, then $z\in S$,
where $x\oplus y$ is the disjoint union of $x$ and $y$;

ii) If $T\subseteq 2^{<\omega}$ is an infinite tree computable in some element of $S$,
then there is $f\in S$ an infinite path through $T$.

\medskip
Scott sets first arose in the study of completions of Peano arithmetic (PA) and models of PA.  Scott \cite{scott} shows that the countable Scott sets are exactly the families of sets ``representable" in a completion of PA.  If $\M$ is a nonstandard model of Peano arithmetic and $a\in\M$, let $$r(a)=\{n\in\omega: \M\models p_n|a\}$$ where $p_0,p_1,\dots$ is an increasing enumeration of the standard primes.  The {\em standard system} of $\M$
 $$SS(\M)=\{r(a):a\in \M\}$$ is a Scott set. 
A longstanding and  vexing problem in the study of models of arithmetic is whether every Scott set arises as the standard system of a model of Peano arithmetic.  The best  result is from Knight and Nadel \cite{kn2}.

\begin{Proposition}\label{kn} If $S$ is a Scott set and $|S|\le \aleph_1$, then there is a model of Peano Arithmetic with standard system $S$.
\end{Proposition}
Thus the Scott set problem has a positive solution if the Continuum \hyphenation {Hypo-thesis} Hypothesis is true, but the question remains open without additional assumptions.   Later in this section, we
sketch a proof of Proposition \ref{kn}.

\medskip
Scott sets also are important when studying recursively saturated structures.
We assume that we are working in a computable language $\L$. We  fix a G\"odel coding of $\L$ and   say that a set of $\L$-formulas is in $S$ if the corresponding set of G\"odel codes is in $S$.

  Let $T$ be a complete $\L$-theory.  The following ideas were introduced in \cite {kn1}, \cite{w} and \cite{mm}.  

\begin{Definition} Let $S\subseteq 2^\omega$.  We say that a model   $\M$ of $T$ is {\em $S$-saturated} if:

i) every type $p\in S_n(\emptyset)$ realized in $\M$ is computable in some element of $S$;

ii) if $p(x,\bar y)\in S_{n+1}(\emptyset)$ is computable in some element of $S$ , $\bar a\in M^n$ and $p(x,\bar a)$ is finitely satisfiable in $\M$, then $p(x,\bar a)$ 
is realized in $\M$.
\end{Definition}

\medskip If a model is $S$-saturated for some $S\subseteq 2^\omega$, then the model is certainly recursively saturated.     

\begin{Proposition}  If $\M\models T$ is recursively saturated, then $\M$ is $S$-saturated for some Scott set $S$.
\end{Proposition}

If the theory $T$ has limited coding power, then we can say little about $S$.  For example, an algebraically closed field of infinite transcendence degree will be $S$-saturated for every Scott set $S$.  On the other hand, the associated Scott set is unique for many natural examples, such as Peano arithmetic, divisible ordered abelian groups, real closed fields, $\Z$-groups
(models of Th$(\Z,+))$ and Presburger arithmetic (models of Th$(\Z,+,<)$).

\begin{Definition} We say that a theory $T$ is {\em effectively perfect} if there is a  tree $(\phi_\sigma:\sigma\in 2^{<\omega})$ of formulas in $n$-free variables computable in $T$ such 
that:

i) $T+\exists \bar v\ \phi_\sigma(\bar v)$ is consistent for all $\sigma$;

ii) if $\sigma\subset \tau$, then $T\models \phi_\tau(\bar v)\rightarrow \phi_\sigma(\bar v)$;

iii) $\phi_{\sigma\^ 0}(\bar v)\land \phi_{\sigma\^ 1}(\bar v)$ is inconsistent with $T$ for all $\sigma$.
\end{Definition}

\begin{Proposition} If $T$ is effectively perfect, then every recursively saturated model of $T$ is $S$-saturated for a unique Scott set $S$.
\end{Proposition}

The theories we will be considering
 are all effectively perfect.
For Peano arithmetic, Presburger arithmetic and $\Z$-groups we can use the formulas $p_n|v$ to find such a tree.
For real closed fields we can use $q<v$  for $q\in \Q$ and for ordered divisible
abelian groups we can use the binary formulas $mv<nw$ for $m,n\in \Z$.

We now sketch a proof of Proposition \ref{kn}.  We first note that for models $\M$ of Peano arithmetic, $\M$ is recursively saturated if and only if $\M$ is $S$-saturated where $S$ is the standard system of $\M$.  (For more details see \cite{kaye}).

\begin{Lemma} If $S$ is a countable Scott set and $T\in S$ is a completion of Peano arithmetic, then there is an $S$-saturated model of $T$.
\end{Lemma} 
\n {\bf Proof Sketch\ } Build $\M$ by a Henkin construction.  At any stage, we will have a finite tuple $\bar a$ and will be committed to $\tp(\bar a)$ the complete type of $\bar a$, where
 $T\subseteq \tp(\bar a)\in S$.  At alternating stages, we either witness an existential quantifier or 
realize a type $p(v,\bar a)\in S$, using the join property of Scott sets to compute
 $p(v,\bar x)\cup  \tp(\bar a)$, and using  the tree property  to find    completions.
\hfill $\Box$

\begin{Lemma} Suppose $S_0\subset S_1$ are countable Scott sets,  $T\in S_0$ is a completion of Peano arithmetic, 
and $\M_0$ and $\M_1$ are countable  recursively saturated models of $T$, where $S_i$ is the standard system of $\M_i$.  Then there is an elementary embedding of $\M_0$ into $\M_1$. 
\end{Lemma}
\n {\bf Proof Sketch} Let $a_0,a_1,\dots$ be a list of the elements of $\M_0$.  Suppose we have a partial elementary map
$(a_0,\dots, a_n)\mapsto (b_0,\dots,b_n)$. If  $\tp(a_{n+1},a_0,\dots,a_n)=p(v,a_0,\dots,a_n)$, there is 
$b\in \M_1$ realizing $p(v,b_0,\dots,b_n)$, and we can extend the embedding.
\hfill $Box$

\medskip We can now prove Proposition \ref{kn}.  Suppose $|S|=\aleph_1$ and $S$ is the union of an $\omega_1$-chain  of countable Scott sets $$S_0\subseteq S_1\subseteq\dots\subseteq S_\alpha\subseteq\dots $$ where $S_\alpha=\bigcup_{\beta<\alpha}S_\beta$ when $\alpha$ is a limit ordinal.  We can build an elementary chain $(\M_\alpha:\alpha<\omega_1)$ where $\M_\alpha$ is recursively saturated with standard system $S_\alpha$. Then $\bigcup_{\alpha<\omega_1} \M_\alpha$ is recursively saturated with standard system  $S$.

While we have nothing new to say about the Scott set problem for Peano arithmetic,
we show that the analgous problem for recursively saturated models has a positive solution in some related algebraic settings.

For divisible ordered abelian groups, this follows easily from an unpublished result of Harnik and Ressayre.  Let $(G,+,<)$ be a divisible ordered abelian group.  Define an equivalence relation on $G$ by $g\equiv h$ if and only if there is a natural number $n$ such that $|g|<n|h|$ and $|h|<n|g|$.  Let $\Gamma=\{|g|/\equiv\ : g\in G\}$,
the set of equivalence classes of positive elements. The ordering of $G$ induces an ordering on $\Gamma$.  Suppose $S$ is a Scott set and $k_S$ is the set of real numbers computable in some element of $S$ (where we identify a real with its cut in the rationals). It is easy to see that $k_S$ is a real closed field.

\begin{Theorem} [Harnik--Ressayre] A divisible ordered abelian group $G$ is $S$-saturated if and only if $\Gamma$ is a dense linear order without endpoints and each equivalence class under $\equiv$ is isomorphic to the ordered additive group of $k_S$.
\end{Theorem}

A complete proof is given in \cite{dkl}.  

\begin{Corollary}\label{odag scott} For any Scott set $S$, there is an $S$-saturated divisible ordered abelian group.
\end{Corollary}
\begin{Proof} Let $G$ be the set of functions $f:\Q\rightarrow k_S$ such that $\{q\in\Q: f(q)\ne 0\}$ is finite.  We add elements of $G$ coordinatewise and order $G$ lexicographically.  
  By the Harnik--Ressayre Theorem, $G$ is $S$-saturated.
\end{Proof}

\section {Real Closed Fields}

In \cite{dks} D'Aquino, Knight and Starchenko show that if $\M$ is a nonstandard model of Peano arithmetic with standard system $S$, then the real closure of the fraction field of $\M$ is an $S$-saturated real closed field.  Thus, it is natural to ask whether we can find an $S$-saturated real closed field for every Scott set $S$.

\begin{Theorem}\label{main rcf} For any Scott set $S$, there is an $S$-saturated real closed field.
\end{Theorem}

 The value group of an $S$-saturated real closed field will be an $S$-saturated divisible ordered abelian group. Thus, Corollary \ref{odag scott} will also follow from   Theorem \ref{main rcf}. 

Theorem \ref{main rcf} is a simple induction using the following Lemma.

\begin{Lemma}  Let $S$ be a Scott set.  Let $K$ be a real closed field such that every type  realized in $K$ is 
in $S$.  Suppose $p(v,\bar w)$ is a set of formulas in $S$, $\bar a\in K$, and $p(v,\bar a)$ is finitely satisfiable in $K$.
Then  we can realize $p(v,\bar a)$ by a (possibly new) element $b$ such that every type realized in $K(b)^\rcl$ is in $S$.
\end{Lemma}

\begin{Proof}   The set of formulas $p(v,\bar a)\cup \tp(\bar a)$ is a consistent partial type in $S$, and, hence, has a completion in $S$. Thus, without loss of generality, we may assume $p(v,\bar a)$ is a complete type.
If $p(v,\bar a)$ is realized in $K$, then there is nothing to do.  If $p(v,\bar a)$ is not realized in $K$
then it determines a cut in the ordering of $\Q(\bar a)^\rcl$ that is not realized in $K$, and, hence, by o-minimality, it
determines a unique type over $K$.  Let $b$ realize $p(v,\bar a)$ and let $\bar c\in K$.  We need to 
show that $\tp(b,\bar c)$ is in $S$.

How do we decide whether $K(b)^{\rm rcl}\models\phi(b,\bar c)$?
By o-minimality, $\phi(v,\bar c)$ defines a finite union of points and intervals 
with endpoints in $\Q(\bar c)^\rcl$.  Since $b\not\in K$, $b$ is neither one of the distinguished points nor an end point of of one of the intervals.
There are $\emptyset$-definable Skolem functions  $f$ and $g$ such 
that:

i) $f(\bar a)<v<g(\bar a)\in \tp(b/\bar a)$;

ii) $f(\bar a)< v<g(\bar a)\rightarrow\phi(v,\bar c)$ or 
$f(\bar a)< v<g(\bar a)\rightarrow\neg\phi(v,\bar c)$;

Given $\tp(b/\bar a)$ and $\tp(\bar a,\bar c)$ we can computably search and find the decomposition of $\phi(v,\bar c)$ and 
 $f$ and $g$ as above.  We can then decide whether $\phi(b,\bar c)$ holds.
Since $S$ is closed under join and Turing reducibility, 
$\tp(b,\bar c)$ is  in $S$.  
\end{Proof}

\medskip  The above argument works for any o-minimal theory $T\in S$.

\medskip  Every real closed field $K$ has a natural valuation for which the valuation ring is $${\cal O}=\{x:|x|<n\hbox{ for some }n\in\Nn\}.$$ If $K$ is recursively saturated, then the value group is a recursively saturated divisible ordered abelian group.  It is natural to ask if every recursively saturated divisible ordered abelian group arises this way.  

D'Aquino, Kuhlmann and Lange \cite{dkl} gave a valuation-theoretic characterization
of recursively saturated real closed fields.  In the following argument we assume familiarity with their results.\footnote{This is essentially our original proof of Theorem 
\ref  {main rcf}.} 

\begin{Proposition} Let $G$ be a recursively saturated divisible ordered abelian group.  There is a recursively saturated real closed field   with value group $G$.
\end{Proposition}
\n {\bf Proof Sketch}  Let $S$ be the Scott set of $G$.  Start with the field $$k_S(t^g:g\in G)^{\rcl}.$$  This is a real closed field
with residue field $k_S$, value group $G$ and all types recursive in $S$.  Given $K$ a real closed field $K$ with value group $G$ and all types recursive in $S$  and suppose we have $\bar a\in K$ and $(f_0,f_1,\dots)$ a sequence of Skolem functions
recursive in $S$ such that $(f_0(\bar a), f_1(\bar a),\dots)$ is pseudo-Cauchy.   If the sequence has no pseudo-limit in $K$ it determines a unique type over $K$.  Adding a realization $b$ does not change the value group or residue field.
As above, every type realized in $K(b)^\rcl$ is in $S$.
We can iterate this construction to build the desired real closed field.
\hfill $\Box$

\section {Presburger Arithmetic}

In \cite{kn2} Knight and Nadel proved that for every Scott set $S$ there is an $S$-saturated $\Z$-group, i.e., an $S$-saturated model of Th$(\Z,+)$.  They asked whether the same is true for       the theory of $(\Z,+,<)$.  This is Presburger arithmetic, which we denote Pr.  We answer this question in the affirmative.

\begin{Theorem}\label{Pr main} For every Scott set $S$, there is an $S$-saturated model of Presburger arithmetic.
\end{Theorem}

We will consider Presburger arithmetic in the language that includes constants for 0 and 1
and unary predicates $P_n(v)$ for $n=2,3,\dots$ that hold if $n$ divides $v$.  We can eliminate quantifiers in this language and the resulting structure is quasi-o-minimal; i.e., any formula $\phi(v,\bar a)$ defines  a finite Boolean combination of $\emptyset$-definable sets
and intervals with endpoints in $\dcl(\bar a)\cup\{\pm\infty\}$. We will use this in the following form. (See, for example \cite{marker} \S 3.1 for quantifier elimination and \cite{bpw} for quasi-o-minimality.)

\begin{Lemma} i) Any formula $\phi(v,\bar a)$ is equivalent over $\tp(\bar a)$ to a Boolean combination of formulas of the form $v\equiv m\mod n$, $v=\alpha$, $v<\beta$ where $\alpha,\beta$ are in the definable closure of $\bar a$ and $m,n\in\N$.

\n ii) $\tp(b,\bar a)$ is determined by:

\begin{itemize} 

\item $\tp(\bar a)$;

\item the sequence $b \mod 2 , b  \mod 3 , b  \mod 4 ,\dots$;

\item the cut of $b$ in the definable closure of $\bar a$.
\end{itemize}
\end{Lemma}  

We obtain Theorem \ref{Pr main} by an iterated construction using the following lemma.

\begin{Lemma} Let $S$ be a Scott set.  Let $G\models \Pr$ such that every type realized in $G$ is computable in  $S$.  Suppose  $\bar a\in G$ and $p(v,\bar w)$
is a complete type  in $S$ such that  $p(v,\bar a)$ is finitely satisfiable.
Then there is $H\supseteq G$ such that $H\models \Pr$, such  that $p(v,\bar a)$ is realized in $G$ and every type realized in $H$ is in $S$.
\end{Lemma}
 \begin{Proof}  If $p(v,\bar a)$ is realized in $G$, then there is nothing to do, so we assume $p(v,\bar a)$ is not realized in $G$. Let  $p^-(v,\bar a)$ be the partial type describing the cut of $v$ over the definable closure of $\bar a$, i.e., $p^-$ consists of all formulas
 of the form $mv < \sum {n_i}a_i$ or $mv>\sum n_i a_i$ that are in $p$ where $m,n_i\in\Z$.
 
 \medskip
\n {\bf Case 1:} Suppose $p^-$ is omitted in $G$.
 
Let $b$ be any realization of $p$, and let $H$ be the definable closure of $G\cup\{b\}$.
It is enough to show that if $\bar c\in G$, then $\tp(b,\bar c)$ is in $S$.    We will show that $\tp(b,\bar c)$ is recursive in $\tp(b,\bar a)$ and $\tp(\bar a,\bar c)$.  Using only  $\tp(b,\bar a)$, we can determine $b  \mod n $ for all $n$.  Thus, we only need to consider formulas of the form $\alpha<v<\beta$
 where $\alpha,\beta\in\dcl(\bar c)$.  Since $p^-(v,\bar a)$ is omitted, we can, as in the case of real closed fields, search to find $\gamma,\delta\in \dcl(\bar a)$ such that $\gamma<b<\delta$ and either
 $\alpha\le\gamma<\delta\le \beta$, $\delta<\alpha$ or $\gamma>\beta$. This can be done recursively in $\tp(b,\bar a) $ and $\tp(\bar a, \bar c)$.
Thus, every type realized in $H$ is in $S$.

\medskip\n {\bf Case 2:}  Suppose  $b\in G$   realizes $p^-$.

Let $\hat b$ be a realization of $p(v,\bar a)$ and let $q_0(v)$ be the divisibility type of $\hat b-b$, i.e., if $\hat b\equiv m \mod n $ and $b\equiv l  \mod n $ then $``v\equiv m-l  \mod n "\in q_0$ for $l,m\in\Z$ and $n>1$.

Let $q(v)\in S_1(G)$ be the unique type containing
\begin{itemize}
\item $q_0(v)$;
\item $n<v$ for all $n\in \Z$;
\item $v<g$ for all $g\in G$ such that $\Z<g$.
\end{itemize}

Let $\epsilon$ realize $q$ and let $H$ be the definable closure of $G\cup\ \{\epsilon\}$.  Suppose $\alpha\in \dcl(\bar a)$ and $b<\alpha$.
Since $G\models$ Pr, we have  $b+n<\alpha$ for all $n\in\Z$. Thus, $\epsilon<\alpha-b$ and $b+\epsilon <\alpha$. Similarly, if $\alpha<b$, then $\alpha<b+\epsilon$ and, thus, $b+\epsilon$ realizes $p^-(v,\bar a)$.  By the choice of $q_0$, $b+\epsilon$ realizes $p(v,\bar a)$.

Suppose $\bar c\in G$.  It suffices to show that $\tp(\epsilon,\bar c)$ is in $S$.
Without loss of generality, we may assume that $\bar c=(c_1,\dots,c_n)$ where all of the $c_i$
are positive infinite and $1,c_1,\dots,c_n$ are linearly independent over $\Q$.
We need to decide the signs of expressions of the form
$$r+s\epsilon+\sum_{i=1}^n t_ic_i$$ where $r,s,t_i\in \Z$.
Such an expression is positive if and only if
\begin{itemize}
\item $\sum t_ic_i>0$, or
\item $\sum t_ic_i=0$ and $q>0$, or
\item $\sum t_ic_i=s=0$ and $r>0$
\end{itemize}
This can be computed using $\tp(\bar c)$.  Thus $\tp(\epsilon,\bar c)$ is recursive in $q_0$
and $\tp(\bar c)$.  Hence, every type realized in $H$ is   in $S$.
\end{Proof}

\end{document}